# A GENERALIZED CHEBYSHEV FINITE DIFFERENCE METHOD FOR HIGHER ORDER BOUNDARY VALUE PROBLEMS

**Soner Aydinlik**

*Department of Mathematics, Faculty of Arts and Sciences, Isik University, 34980, Istanbul, Turkey*

**Ahmet Kiris**[*]

*Department of Mathematical Engineering, Faculty of Arts and Sciences, Istanbul Technical University, 34469, Istanbul, Turkey*

**Abstract**

A general formula is presented for any order derivative of Chebyshev polynomials instead of the existing recursive relationship. Hence, the Chebyshev finite difference method is made applicable not only to second order problems but also to higher order boundary value problems. The generalized method is applied to a variety of higher order boundary value problems and it is seen that the obtained results are more accurate than the other numerical methods in absolute error.

**Keywords:** Chebyshev finite difference method, numerical solution, higher order boundary value problems.

## 1. Introduction

Many engineering problems including dynamical system problems such as harmonic oscillator, elasticity problems such as wave propagation and heat convection, and calculus of variations problems have been modelled with second-order boundary value problems. However, when mathematical constraints are slightly stretched in order to provide compatibility with physical realities or when the factors that affect the problem are analyzed concomitantly, higher order boundary value problems arise.

There exist some methods such as finite difference and the shooting method to solve second-order boundary value problems, unlikely, more effective methods are required to solve higher order boundary value problems; Aouadi used the Chebyshev finite difference method to solve the third-order boundary value problem arising in the modelling of mass transfer when the material was considered as a micropolar material [1, 2] instead of a classical elastic material [3]. Fifth-order boundary value problems arise in the mathematical modelling of viscoelastic fluids [4, 5]. To solve these problems, Khan used the finite difference method [6], Çağlar et al. used the sixth-degree B-spline and collocation methods [7], and Wazwaz used the Adomian method [8]. In order to solve six-degree problems occurring in the modelling of many astrophysics problems, septic spline [9], Legendre-Galerkin [10] and Daftardar-Jafari methods [11] were used. To solve eighth-order boundary value problems that arise in torsional vibration of uniform beams, the Adomian decomposition method [12], the differential quadrature method [13], the homotopy perturbation method [14], and the spline method [15] were used. When fluid layers are subject to rotation, heat convection is modelled with $10^{th}$ order differential equations, if magnetic effects are also included [16, 17]. Ullan, Khan, and Rahim used a new iterative technique [18] to solve such problems.

[*] Corresponding Author; kiris@itu.edu.tr

The Chebyshev finite difference method is more advantageous in solution of higher order problems than the methods mentioned above. Of these methods, for example the Adomian method requires calculation of Adomian polynomials and the homotopy methods require many conditions in addition to finding the appropriate parameters [19]. Because of their orthogonality property, the Chebyshev polynomials form a complete orthogonal set on a space of continuous functions and since recursive relations can be obtained easily, especially derivatives can be calculated recursively at any order. Furthermore, not only values at specific points of the solution range -as in many numerical methods-, but an approximation polynomial valid throughout all the interval is obtained. The greatest advantage of Chebyshev polynomials compared to the other polynomial approximations of the same order is that they are the most convenient polynomials having the lowest maximum error in the given range [20].

Solution of second-order initial or boundary value problems with the Chebyshev finite difference method is widely used in the literature; El-Kady and Elbarbary obtained a general formula for the derivatives up to second order of the Chebyshev approximation polynomial instead of the recursive relation and solved second-order boundary value problems with the help of this formula [21]. Saadatmandi and Farsangi solved the second-order nonlinear system and Saadatmandi and Deghan solved some calculus of variation problems with the Chebyshev finite difference method; however, in these studies, only second order problems could be solved [22, 23]. Aouadi obtained another formula containing successive sums for the third order derivative to analyze micropolar flow and mass transfer from a surface stretched with heat [3]. However, successive sums increase complexity as the order increases. This is the main reason of why Chebyshev finite difference method could not be used for higher order differential equations.

In this study, first, the derivatives of the Chebyshev polynomial of any order are obtained without any recursive relation and then by the help of these derivatives, a general formula is presented for the Chebyshev approximation polynomial. With this formula, the Chebyshev finite difference method became applicable not only to first and second order problems, but also to initial or boundary value problems of any order. The generalized Chebyshev finite difference method is applied to a variety of higher order boundary value problems given in the literature and the obtained results are more accurate than the other numerical methods in absolute error.

## 2. Chebyshev Polynomials

Chebyshev polynomials of the first kind are defined as [20]:

$$T_n(x) = \cos(n\theta), \quad x = \cos\theta, \quad \theta \in [0,\pi], \quad \forall n \geq 0, \quad x \in [-1,1]. \tag{2.1}$$

Because of their trigonometric properties, the recursive relation

$$T_{n+1}(x) = 2x\,T_n(x) - T_{n-1}(x) \tag{2.2}$$

can be easily observed. Chebyshev polynomials are orthogonal in the interval $[-1,1]$ respect to the weight function, $w(x) = \dfrac{1}{\sqrt{1-x^2}}$. When $n \geq 1$ in the $x \in [-1,1]$ they have $n$ roots at the points

$$x_k = \cos\left(\frac{2k-1}{2n}\pi\right), \quad k = 1, 2, \ldots, n \tag{2.3}$$

and at the points

$$\tilde{x}_k = \cos\left(\frac{k\pi}{n}\right), \quad k = 0, 1, 2, \ldots, n \tag{2.4}$$

they have $(n+1)$ extrema. An important inequality, which is known as the economization property of the Chebyshev polynomials is given by the theorem given below.

**Theorem:** Let $\tilde{\prod}_n$ are the sets of all monic polynomials of the order $n$ and the monic Chebyshev polynomials are defined by $\tilde{T}_n(x) = \dfrac{T_n(x)}{2^{n-1}}$. Then the inequality

$$\frac{1}{2^{n-1}} = \|\tilde{T}_n(x)\|_\infty \leq \|P_n(x)\|_\infty \tag{2.5}$$

holds for $\forall P_n(x) \in \tilde{\Pi}_n$, here the equality is satisfied only when $P_n(x) \equiv \tilde{T}_n(x)$. The proof can be find in [20]. The economization property states that when approaching a function $f(x)$ with polynomials, in order to minimize the maximum error in the given range, the Gauss-Lobatto points, i.e. the roots of the polynomial $\tilde{T}_{n+1}(x)$ have to be taken as the node points in the interpolation. This is expressed with

$$\max_{x \in [-1,1]} |f(x) - P_n(x)| \leq \frac{1}{2^n (n+1)!} \max_{x \in [-1,1]} |f^{(n+1)}(x)|. \tag{2.6}$$

## 3. Chebyshev Finite Difference Method

Clenshaw and Curtis defined the solution for a given initial or boundary value problem as a series of Chebyshev polynomials [24]

$$y(x) = \sum_{n=0}^{N} {}''a_n T_n(x) \tag{3.1}$$

where, the superscript ($''$) in the sum symbol means that half of the first and the last terms have to be taken and $N$ shows the order of the approximation polynomial. Using the orthogonality property of Chebyshev polynomials and the Lagrange interpolation, the unknown $a_n$ coefficients are found as

$$a_n = \frac{2}{N} \sum_{j=1}^{N} y(x_j) T_n(x_j). \tag{3.2}$$

The value of the $m^{th}$ order derivative of $y(x)$ in (3.1) at the points $x_k$ is given by

$$y^{(m)}(x_k) = \sum_{j=0}^{N} d_{k,j}^{(m)} y(x_j) \tag{3.3}$$

where,

$$d_{k,j}^{(m)} = \frac{2}{N} \sum_{n=0}^{N} {}''T_n(x_j) T_n^{(m)}(x_k) \tag{3.4}$$

and only the derivatives of Chebyshev polynomials, $T_n^{(m)}(x)$ are unknowns. Elbarbary and El-Kady showed that the first two of the coefficients, $d_{k,j}^{(m)}$ are given as [21]

$$d_{k,j}^{(1)} = \frac{4\theta_j}{N} \sum_{n=0}^{N} \sum_{\substack{l=0 \\ (n+l)\,odd}}^{n-1} \frac{n\theta_n}{c_l} T_n(x_j) T_l(x_k), \quad k, j = 0,1,...,N, \tag{3.5}$$

$$d_{k,j}^{(2)} = \frac{2\theta_j}{N} \sum_{n=0}^{N} \sum_{\substack{l=0 \\ (n+l)\,even}}^{n-2} \frac{n(n^2-l^2)\theta_n}{c_l} T_n(x_j) T_l(x_k), \quad k, j = 0,1,...,N \tag{3.6}$$

where

$$\left.\begin{array}{l} \theta_0 = \theta_n = 1/2, \; \theta_j = 1 \\ c_0 = 2, \; c_i = 1 \end{array}\right\} \quad j = 1,2,...,N-1, \quad i \geq 1. \tag{3.7}$$

On the other hand, Aouadi gives the coefficients for the third order as [3]

$$d_{k,j}^{(3)} = \frac{4\theta_j}{L} \sum_{n=0}^{L} \sum_{\substack{l=0 \\ (n+l)\,even}}^{n-2} \sum_{\substack{i=0 \\ (i+l)\,odd}}^{l-1} \frac{\theta_n n l}{c_l c_i}(n^2-l^2)^2 T_n(x_j) T_l(x_k), \quad (n+l)\;even,\;(i+l)\;odd. \tag{3.8}$$

As it can be seen in (3.8), a new sum symbol comes for third order. Emerging a new sum symbol for every order in (3.8) and not knowing $T_n^{(m)}(x)$ in (3.4) make the calculation of higher order coefficients impossible in a similar way. In the presented study, in order to make the coefficients $d_{k,j}^{(m)}$ computable for each order, the relation between Chebyshev polynomials

$$T_n(x) = \begin{cases} \dfrac{T'_{n+1}(x)}{n+1}, & n = 0 \\[2mm] \dfrac{T'_{n+1}(x)}{2(n+1)}, & n = 1 \\[2mm] \dfrac{T'_{n+1}(x)}{2(n+1)} - \dfrac{T'_{n-1}(x)}{2(n-1)}, & n > 1 \end{cases} \tag{3.9}$$

are used, and then the recursive relation is eliminated by taking successive derivatives of (3.9). Consequently, the derivatives of the Chebyshev polynomials are found as

$$T_n^{(3)}(x) = \sum_{\substack{l=0 \\ (n+l)\,odd}}^{n-3} \frac{1}{4c_l} n(n^2-(l+1)^2)(n^2-(l-1)^2) T_l(x),$$

$$T_n^{(4)}(x) = \sum_{\substack{l=0 \\ (n+l)\,even}}^{n-4} \frac{1}{24c_l} n(n^2-(l-2)^2)(n^2-l^2)(n^2-(l+2)^2) T_l(x), \tag{3.10}$$

$$T_n^{(5)}(x) = \sum_{\substack{l=0 \\ (n+l)\,odd}}^{n-5} \frac{1}{192c_l} n(n^2-(l-3)^2)(n^2-(l-1)^2)(n^2-(l+1)^2)(n^2-(l+3)^2) T_l(x)$$

and following a similar process as above, finally the derivative of any order is given as [25]

$$T_n^{(m)}(x) = \sum_{\substack{l=0 \\ (n+l+m)\text{ even}}}^{n-m} \prod_{\substack{i=2-m \\ m>1}}^{m-2} {}^*\left(n^2 - (l+i)^2\right) \frac{n}{c_l} \frac{1}{(m-1)!\, 2^{(m-2)}} T_l(x), \quad m \geq 1. \tag{3.11}$$

From this general formula, the value of the $m^{th}$ order derivative at the points $x_k$, $T_n^m(x_k)$ is calculated and substituting this into (3.4) gives the coefficients $d_{k,j}^{(m)}$ as

$$d_{k,j}^{(m)} = \frac{\theta_j}{N} \sum_{n=0}^{N} \sum_{\substack{l=0 \\ (n+l+m)\text{ even}}}^{n-m} \prod_{i=2-m}^{m-2} {}^*\left(n^2 - (l+i)^2\right) \frac{\theta_n n}{c_l} \frac{1}{(m-1)!\, 2^{(m-3)}} T_n(x_j) T_l(x_k), \quad m>1, \quad k,j=0,1,\ldots,N \tag{3.12}$$

The symbol (*) in the multiplication shows that the multiplication index increases two by two in both (3.11) and (3.12). Thus, the $m^{th}$ order derivative of the approximation polynomial may be easily calculated in (3.3) by the help of (3.12).

Finally, writing the given initial or boundary value problem in terms of approximation polynomial and its derivatives, the problem is transformed to a linear or nonlinear algebraic equation. Evaluating this equation and the initial or boundary conditions at the Gauss-Lobatto points, i.e., $(x_j,\ j=0,1,\ldots,N)$, a system consisting of $(N+1)$ equations is obtained. This system can be solved using any appropriate method and the approximation polynomial's $(N+1)$ unknowns, which constitute the solution to the original differential equation, are determined.

The Chebyshev finite difference method is made applicable not only for the first two order but also to higher order initial or boundary value problems by the generalized formulas (3.11) and (3.12). In the following section, the method is applied to some higher order boundary value problems and compared to the methods used commonly in the literature. It is shown that the error of the presented method is much lower than those of other methods.

## 4. Numerical Examples

In this section, five nonlinear problems are solved by generalized Chebyshev finite difference method mentioned above.

**Example 1:** Consider the nonlinear boundary value problem

$$y^{(4)}(x) = \sin x + \sin^2 x - (y''(x))^2 \quad 0 < x < 1, \tag{4.1}$$

with the given boundary conditions

$$y(0) = 0, \quad y'(0) = 1,$$
$$y(1) = \sin 1, \quad y'(1) = \cos 1. \tag{4.2}$$

The exact solution for the above problem is $y(x) = \sin x$. In order to solve this problem with using the Chebyshev finite difference method, the interval of $[0,1]$ should be transferred to the interval of $[-1,1]$ by $t = 2x - 1$. The transferred differential equation and the boundary conditions are

$$2^4 y^{(4)}(t) = \sin\left(\frac{(t+1)}{2}\right) + \sin^2\left(\frac{(t+1)}{2}\right) - 2^4 \left(y''(t)\right)^2, \quad -1 < t < 1, \tag{4.3}$$

$$y(-1) = 0, \quad y'(-1) = \frac{1}{2},$$
$$y(1) = \sin 1, \quad y'(1) = \frac{\cos 1}{2}. \tag{4.4}$$

(4.3) and (4.4) are converted to the following form with the help of (3.3)

$$2^4 \sum_{j=0}^{N} d_{k,j}^4 y(t_j) = \sin\left(\frac{t_k+1}{2}\right) + \sin^2\left(\frac{t_k+1}{2}\right) - 2^4 \left(\sum_{j=0}^{N} d_{k,j}^2 y(t_j)\right)^2, \quad k = 1, 2, ..., N-1 \tag{4.5}$$

$$y(t_N) = 0, \quad \sum_{j=0}^{N} d_{N,j}^1 y(t_j) = \frac{1}{2},$$
$$y(t_0) = \sin 1, \quad \sum_{j=0}^{N} d_{0,j}^1 y(t_j) = \frac{\cos 1}{2}. \tag{4.6}$$

Equations (4.5) and (4.6) gives $N-1$ nonlinear algebraic equations containing $N-1$ unknown coefficients, $y(t_1), y(t_2), ..., y(t_{N-1})$ and can be solved by any appropriate numerical root finding method. Here, the obtained values of $y(t)$ from the system (4.5) and (4.6) constitutes the solution of the (4.3) in virtue of (3.1) and (3.2). Using the inverse transformation $x = \frac{t+1}{2}$, the solution of the given boundary value problem, $y(x)$ is obtained as

$$\begin{aligned}
y(x) = & 1.266764471097303 \times 10^{-13} + 0.9999999999629425\, x \\
& + 2.231311690970017 \times 10^{-9}\, x^2 - 0.16666671897218904\, x^3 \\
& + 6.343300914361805 \times 10^{-7}\, x^4 + 0.008328746101867068\, x^5 \\
& + 0.00002142878559295363\, x^6 - 0.00026621823852630253\, x^7 \\
& + 0.00014926375804082324\, x^8 - 0.00022833703094095526\, x^9 \\
& + 0.0002506193902393404\, x^{10} - 0.00018630278230245625\, x^{11} \\
& + 0.00009032584992902619\, x^{12} - 0.000025715279792036328\, x^{13} \\
& + 0.00000325670199734824\, x^{14}.
\end{aligned} \quad (4.7)$$

The approximation polynomial, $y(x)$ is not given for the other examples for the sake of brevity.

The comparison of the results with the variational iteration method [26] and Chebyshev finite difference method is given in Table 3.11. As it can be seen from the Table, the absolute error between the analytical solution and the result obtained by CFDM is less than the absolute error of the given method in [26]. Furthermore, the default precision in Mathematica, which is used here is 16, if the precision is setup to 100, the results is much more convincing.

**Table 3.11:** The comparison of example 1, the analytical solution is $\sin x$.

| $x$ | Absolute Error by VIM | Absolute Error by CFDM |
|---|---|---|
| 0.0 | $9.5923E-14$ | $1.2669E-13$ |
| 0.1 | $7.7856E-08$ | $7.7729E-14$ |
| 0.2 | $2.7231E-07$ | $6.2561E-14$ |
| 0.3 | $5.2489E-07$ | $2.5512E-13$ |
| 0.4 | $7.7730E-07$ | $1.3716E-13$ |
| 0.5 | $9.7145E-07$ | $1.5559E-13$ |
| 0.6 | $1.0502E-06$ | $1.6986E-13$ |
| 0.7 | $9.6286E-07$ | $3.5693E-13$ |
| 0.8 | $6.8407E-07$ | $4.5319E-13$ |
| 0.9 | $2.7069E-07$ | $5.8797E-13$ |
| 1.0 | $1.5676E-13$ | $4.916E-13$ |

**Example 2:** Consider the nonlinear boundary value problem

$$y^{(6)}(x) - 20\, e^{-36 y(x)} = -40\,(1+x)^{-6}, \quad 0 < x < 1, \tag{4.8}$$

with the given boundary conditions

$$y(0) = 0, \quad y'(0) = \frac{1}{6}, \quad y''(0) = \frac{-1}{6},$$
$$y(1) = \frac{1}{6}\ln 2, \quad y'(1) = \frac{1}{12}, \quad y''(1) = \frac{-1}{24}. \tag{4.9}$$

The exact solution for the above problem is $y(x) = \frac{1}{6}\ln(1+x)$. Similar to Example 1, if the linear transformation $t \to 2x - 1$ is used, the following equations are obtained

$$2^6\, y^{(6)}(t) - 20\, e^{-36 y(t)} = -40\left(\frac{t+3}{2}\right)^{-6}, \quad -1 < t < 1, \tag{4.10}$$

$$y(-1) = 0, \quad y'(-1) = \frac{1}{12}, \quad y''(-1) = \frac{-1}{24},$$
$$y(1) = \frac{1}{6}\ln 2, \quad y'(1) = \frac{1}{24}, \quad y''(1) = \frac{-1}{96}. \tag{4.11}$$

(4.10) and (4.11) are reduced to the system of algebraic equations with the help of (3.3) as

$$2^6 \sum_{j=0}^{N} d^6_{k,j}\, y(t_j) - e^{-36 y(t_k)} = -40\left(\frac{t_k + 3}{2}\right)^{-6}, \quad k = 1, 2, \ldots, N-1 \tag{4.12}$$

$$y(t_N) = 0, \quad \sum_{j=0}^{N} d^1_{N,j}\, y(t_j) = \frac{1}{12}, \quad \sum_{j=0}^{N} d^2_{N,j}\, y(t_j) = \frac{-1}{24},$$
$$y(t_0) = \frac{1}{6}\ln 2, \quad \sum_{j=0}^{N} d^1_{0,j}\, y(t_j) = \frac{1}{24}, \quad \sum_{j=0}^{N} d^2_{0,j}\, y(t_j) = \frac{-1}{96}. \tag{4.13}$$

Using the solution of the system (4.12) and (4.13) in (3.1) and (3.2) after the inverse transformation gives the approximation Chebyshev polynomial which is the solution of the given problem. The comparison of the results with the quintic B-spline collocation [27] and Chebyshev finite difference methods is given in Table 3.12.

**Table 3.12:** The comparison of example 2, the analytical solution is $\frac{1}{6}\ln(1+x)$.

| $x$ | Absolute Error by QBSCM | Absolute Error by CFDM |
| --- | --- | --- |
| 0.1 | $1.303852E - 07$ | $7.235878E - 14$ |

| 0.2 | 5.215406E−07 | 2.065990E−13 |
| 0.3 | 9.723008E−07 | 2.936120E−13 |
| 0.4 | 1.329929E−06 | 2.937860E−13 |
| 0.5 | 1.259148E−06 | 2.183110E−13 |
| 0.6 | 8.419156E−07 | 1.075670E−13 |
| 0.7 | 4.023314E−07 | 1.088020E−14 |
| 0.8 | 8.195639E−08 | 3.384792E−14 |
| 0.9 | 1.713634E−07 | 2.166322E−14 |

**Example 3:** Consider the nonlinear boundary value problem

$$y^{(6)}(x) = e^x y(x)^2, \quad 0 < x < 1, \tag{4.14}$$

with the given boundary conditions

$$\begin{aligned} y(0) &= 1, \quad y'(0) = -1, \quad y''(0) = 1, \\ y(1) &= e^{-1}, \quad y'(1) = -e^{-1}, \quad y''(1) = e^{-1}. \end{aligned} \tag{4.15}$$

The exact solution for the above problem is $y(x) = e^{-x}$. With the same transformation above, the following equations are obtained

$$2^6 y^{(6)}(t) = e^{\frac{t+1}{2}} y(t)^2, \quad -1 < t < 1 \tag{4.16}$$

$$\begin{aligned} y(-1) &= 1, \quad y'(-1) = \frac{-1}{2}, \quad y''(-1) = \frac{1}{4}, \\ y(1) &= e^{-1}, \quad y'(1) = \frac{-e^{-1}}{2}, \quad y''(1) = \frac{e^{-1}}{4}. \end{aligned} \tag{4.17}$$

Similarly, (4.16) and (4.17) are transformed to

$$2^6 \sum_{j=0}^{N} d_{k,j}^6 y(t_j) = e^{\frac{t_k+1}{2}} y(t_k)^2, \quad k = 1, 2, \ldots, N-1 \tag{4.18}$$

$$\begin{aligned} y(t_N) &= 1, \quad \sum_{j=0}^{N} d_{N,j}^1 y(t_j) = \frac{-1}{2}, \quad \sum_{j=0}^{N} d_{N,j}^2 y(t_j) = \frac{1}{4}, \\ y(t_0) &= e^{-1}, \quad \sum_{j=0}^{N} d_{0,j}^1 y(t_j) = \frac{-e^{-1}}{2}, \quad \sum_{j=0}^{N} d_{0,j}^1 y(t_j) = \frac{e^{-1}}{4}. \end{aligned} \tag{4.19}$$

Same process is done to find the approximation polynomial. The comparison of the results with the Adomian [28], homotopy perturbation [29], variational iteration [30], new iterative, Daftardar Jafari [31] and Chebyshev finite difference method is given in Table 3.13.

**Table 3.13**: The comparison of example 3, the analytical solution is $e^{-x}$.

| $x$ | Absolute Error[ADM] | Absolute Error [HPM] | Absolute Error [VIM] | Absolute Error [ITM] | Absolute Error [DJM] | Absolute Error [CFDM] |
|---|---|---|---|---|---|---|
| 0.1 | $2.4E-7$ | $2.4E-7$ | $2.4E-7$ | $2.4E-7$ | $3.1E-14$ | $3.3E-15$ |
| 0.2 | $1.4E-6$ | $1.4E-6$ | $1.4E-6$ | $1.4E-6$ | $1.9E-13$ | $2.0E-14$ |
| 0.3 | $3.3E-6$ | $3.3E-6$ | $3.3E-6$ | $3.3E-6$ | $4.8E-13$ | $5.4E-14$ |
| 0.4 | $5.2E-6$ | $5.2E-6$ | $5.2E-6$ | $5.2E-6$ | $8.0E-13$ | $9.1E-14$ |
| 0.5 | $6.2E-6$ | $6.2E-6$ | $6.2E-6$ | $6.2E-6$ | $1.0E-12$ | $1.2E-13$ |
| 0.6 | $5.8E-6$ | $5.8E-6$ | $5.8E-6$ | $5.8E-6$ | $1.0E-12$ | $1.2E-13$ |
| 0.7 | $4.1E-6$ | $4.1E-6$ | $4.1E-6$ | $4.1E-6$ | $8.1E-13$ | $1.0E-13$ |
| 0.8 | $1.9E-6$ | $1.9E-6$ | $1.9E-6$ | $1.9E-6$ | $4.3E-13$ | $6.3E-14$ |
| 0.9 | $3.6E-7$ | $3.6E-7$ | $3.6E-7$ | $3.6E-7$ | $9.2E-13$ | $1.6E-14$ |
| 1.0 | $5.0E-10$ | $5.0E-10$ | $5.0E-10$ | $5.0E-10$ | $5.6E-14$ | $0$ |

**Example 4:** Consider the problem

$$y^{(7)}(x) + y^{(4)}(x) - e^{y(x)} y(x) = e^x((12 - 4x + (x-1) e^{(-e^x(x-1)\cos x)})\cos x \\ - 8(5+x)\sin x), \quad 0 < x < 1, \tag{4.20}$$

with the conditions given

$$y(0) = 1, \quad y'(0) = 0, \quad y''(0) = -2, \quad y'''(0) = -2, \\ y(1) = 0, \quad y'(1) = -e\cos 1, \quad y''(1) = -2e\cos 1 + 2\sin 1. \tag{4.21}$$

The analytical solution for the problem is $y(x) = e^x(1-x)\cos(x)$. The same procedure gives (4.20) and (4.21) as

$$2^7 y^{(7)}(t) + 2^4 y^{(4)}(t) - e^{y(t)} y(t) = e^{\frac{t+1}{2}}((10 - 2t + \frac{t-1}{2} e^{(-\frac{t-1}{2} e^{\frac{t+1}{2}} \cos(\frac{t+1}{2}))})\cos(\frac{t+1}{2}) \\ - 8\frac{t+11}{2} \sin(\frac{t+1}{2})), \quad -1 < t < 1, \tag{4.22}$$

$$y(-1) = 1, \quad y'(-1) = 0, \quad y''(-1) = -\frac{1}{2}, \quad y'''(-1) = -\frac{1}{4},$$
$$y(1) = 0, \quad y'(1) = -\frac{e\cos 1}{2}, \quad y''(1) = \frac{-e\cos 1 + \sin 1}{2}.$$
(4.23)

Using (3.3) reduces (4.22) and (4.23) to

$$2^7 \sum_{j=0}^{N} d_{k,j}^7 y(t_j) + 2^4 \sum_{j=0}^{N} d_{k,j}^4 y(t_j) - e^{y(t_k)} y(t_k) = e^{\frac{t_k+1}{2}} ((10 - 2t_k + \frac{t_k-1}{2} e^{(-e^{\frac{t_k+1}{2}} \frac{t_k-1}{2} \cos(\frac{t_k+1}{2}))})$$
$$\times \cos(\frac{t_k+1}{2}) - 8\frac{t_k+11}{2}\sin(\frac{t_k+1}{2})), \qquad k = 1, 2, ..., N-1$$
(4.24)

$$y(t_N) = 1, \quad \sum_{j=0}^{N} d_{N,j}^1 y(t_j) = 0, \quad \sum_{j=0}^{N} d_{N,j}^2 y(t_j) = -\frac{1}{2}, \quad \sum_{j=0}^{N} d_{N,j}^3 y(t_j) = -\frac{1}{4},$$
$$y(t_0) = 0, \quad \sum_{j=0}^{N} d_{0,j}^1 y(t_j) = -\frac{e\cos 1}{2}, \quad \sum_{j=0}^{N} d_{0,j}^2 y(t_j) = \frac{-e\cos 1 + \sin 1}{2}.$$
(4.25)

The comparison of the reproducing kernel space [32] and Chebyshev finite difference method is given in Table 3.14.

**Table 3.14:** The comparison of example 4, the analytical solution is $e^x(1-x)\cos(x)$.

| $x$ | Absolute Error by RKS with $n=30$ | Absolute Error by RKS with $n=50$ | Absolute Error by CFDM |
|---|---|---|---|
| 0.0 | $6.4278E-11$ | $6.4113E-11$ | $1.04949E-12$ |
| 0.125 | $4.7378E-10$ | $1.4645\ E-10$ | $2.06599E-12$ |
| 0.250 | $5.2047E-09$ | $1.9111E-09$ | $1.10497E-11$ |
| 0.375 | $1.5281E-08$ | $5.6158E-09$ | $2.36455E-11$ |
| 0.500 | $2.4509E-08$ | $8.8518E-09$ | $2.11277E-11$ |
| 0.625 | $2.5265E-08$ | $9.1373E-09$ | $1.13554E-12$ |
| 0.750 | $1.5563E-08$ | $5.6666E-09$ | $1.15951E-11$ |
| 0.875 | $3.2941E-09$ | $1.0112E-09$ | $4.8044E-12$ |
| 1.0 | $5.6254E-11$ | $5.6239E-11$ | $5.05727E-13$ |

**Example 5:** Consider the problem

$$y^{(10)}(x) = \frac{14175}{4}(x + y(x) + 1)^{11}, \qquad 0 < x < 1$$
(4.26)

with the conditions

$$y(0) = 0, \quad y'(0) = \frac{-1}{2}, \quad y''(0) = \frac{1}{2}, \quad y^{(3)}(0) = \frac{3}{4}, \quad y^{(4)}(0) = \frac{3}{2}, \quad (4.27)$$

$$y(1) = 0, \quad y'(1) = 1, \quad y''(1) = 4, \quad y^{(3)}(1) = 12, \quad y^{(4)}(1) = 48.$$

The analytical solution is $y(x) = \frac{2}{2-x} - x - 1$. With the same procedure one can get the following equations

$$2^{10} y^{(10)}(t) = \frac{14175}{4}\left(\frac{t+3+2y(t)}{2}\right)^{11}, \quad -1 < t < 1, \quad (4.28)$$

$$y(-1) = 0, \quad y'(-1) = \frac{-1}{4}, \quad y''(-1) = \frac{1}{8}, \quad y^{(3)}(-1) = \frac{3}{32}, \quad y^{(4)}(-1) = \frac{3}{32},$$

$$y(1) = 0, \quad y'(1) = \frac{1}{2}, \quad y''(1) = 1, \quad y^{(3)}(1) = \frac{3}{2}, \quad y^{(4)}(1) = 3. \quad (4.29)$$

These equations are transformed to the following form

$$2^{10} \sum_{j=0}^{N} d_{k,j}^{10} y(t_j) = \frac{14175}{4}\left(\frac{t_k + 3 + 2y(t_k)}{2}\right)^{11}, \quad k = 1, 2, \ldots, N-1 \quad (4.30)$$

$$y(t_N) = 0, \quad \sum_{j=0}^{N} d_{N,j}^1 y(t_j) = \frac{-1}{4}, \quad \sum_{j=0}^{N} d_{N,j}^2 y(t_j) = \frac{1}{8}, \quad \sum_{j=0}^{N} d_{N,j}^3 y(t_j) = \frac{3}{32}, \quad \sum_{j=0}^{N} d_{N,j}^4 y(t_j) = \frac{3}{32},$$

$$y(t_0) = 0, \quad \sum_{j=0}^{N} d_{0,j}^1 y(t_j) = \frac{1}{2}, \quad \sum_{j=0}^{N} d_{0,j}^2 y(t_j) = 1, \quad \sum_{j=0}^{N} d_{0,j}^3 y(t_j) = \frac{3}{2}, \quad \sum_{j=0}^{N} d_{0,j}^4 y(t_j) = 3. \quad (4.31)$$

When the inverse transformation is applied to the numerical solution of the system (4.30) and (4.31), and then substituting this into (3.1) and (3.2) gives the approximation polynomial. The comparison of the result with the quintic B-spline collocation method [33] and homotopy analysis method [34] is given in Table 3.15.

**Table 3.15:** The comparison of example 5, the analytical solution is $\frac{2}{2-x} - x - 1$.

| $x$ | Absolute Error by QBSCM | Absolute Error by HAM | Absolute Error by CFDM |
|---|---|---|---|
| 0.1 | $1.322478E-06$ | $3.95413E-11$ | $1.335263E-11$ |
| 0.2 | $4.231930E-06$ | $7.33317E-10$ | $2.468730E-10$ |
| 0.3 | $1.676381E-05$ | $7.33317E-09$ | $8.750325E-10$ |
| 0.4 | $4.245341E-05$ | $6.06524E-09$ | $1.393045E-09$ |
| 0.5 | $6.662799E-05$ | $7.74775E-09$ | $1.202488E-09$ |
| 0.6 | $6.940961E-05$ | $6.56402E-09$ | $5.293163E-10$ |
| 0.7 | $4.750490E-05$ | $3.48667E-09$ | $5.506837E-11$ |
| 0.8 | $1.643598E-05$ | $9.23198E-10$ | $2.995830E-11$ |
| 0.9 | $2.607703E-07$ | $5.33521E-11$ | $4.110150E-12$ |

## 5. Conclusion

A general non-recursive formula is presented for the derivatives of the Chebyshev polynomials of any order instead of existing recursive relationship of the Chebyshev polynomials. By using these derivatives, a general formula is obtained for the Chebyshev approximation polynomial. Thus, the Chebyshev finite difference method is made applicable to not only the first and second order problems as in the previous studies, but also higher order initial or boundary value problems. The obtained results show that this generalization of the Chebyshev method can solve boundary value problems efficiently, and the better accuracy is observed in comparison with the presented method and existing techniques. For future prospects, the method is expected to be applied to partial differential equations for higher order problems.